\documentclass{commat}

\usepackage{graphicx}

\title{%
   On the Diophantine equation $B_{n_{1}}+B_{n_{2}}=2^{a_{1}}+2^{a_{2}}+2^{a_{3}}$
    }

\author{%
   Kisan Bhoi and Prasanta Kumar Ray
    }

\affiliation{
    \address{Kisan Bhoi --
    Sambalpur University, Jyoti Vihar, Burla, India
        }
    \email{%
   kisanbhoi.95@suniv.ac.in
    }
    \address{Prasanta Kumar Ray --
   Sambalpur University, Jyoti Vihar, Burla, India
        }
    \email{%
   prasantamath@suniv.ac.in
    }
    }

\abstract{%
   In this study we find all solutions of the Diophantine equation
   \[
   		B_{n_{1}}+B_{n_{2}}=2^{a_{1}}+2^{a_{2}}+2^{a_{3}}
   	\]
   	in positive integer variables $(n_{1},n_{2},a_{1},a_{2},a_{3}),$ where $B_{n}$ denotes the $n$-th balancing number.
    }



\VOLUME{31}
\YEAR{2023}
\NUMBER{1}
\firstpage{375}
\DOI{https://doi.org/10.46298/cm.10476}

\begin{paper} 

\section{Introduction}
Balancing sequence $\{B_{n}\}_{n\geq1}$ is originated from a simple Diophantine equation 
\[
	1+2+\ldots+(n-1)=(n+1)+(n+2)+\cdots+(n+r)
\]
introduced by Behera and Panda \cite{BEHERA}. Here, $r$ is called a balancer corresponding to a balancing number $n.$ The balancing sequence satisfies the binary recurrence
\begin{align*}
B_{n+1}=6B_{n}-B_{n-1},\quad  n\geq 1
\end{align*}
with seeds $B_{0}=0$ and $B_{1}=1.$ The Binet's formula for $\{B_{n}\}_{n\geq1}$ is given by
\begin{align*}
B_{n}=\frac{\alpha^{n}-\beta^{n}}{4\sqrt{2}},
\end{align*}
where $\alpha=3+\sqrt{8}$ and $\beta=3-\sqrt{8}$ are the zeros of the polynomial $f(x)=x^{2}-6x+1.$
Clearly, $\beta^{-1}=\alpha.$ It can be easily seen that
\begin{align}\label{eq:1**}
\alpha^{n-1}<B_{n}<\alpha^{n},~~~\text{for}~~~n>1.
\end{align}
Diophantine equations involving powers and binary recurrence sequences have been extensively studied by many researchers in recent past. For example, Bravo and Luca \cite{BRAVO} found all solutions of the equation $F_{n}+F_{m}=2^{a},$ where $F_{n}$ is the $n$-th Fibonacci number. Later, Bravo and Bravo \cite{EFBRAVO} extended this work and found all positive integer solutions of the Diophantine equation $F_{n}+F_{m}+F_{l}=2^{a}.$ In \cite{KESKIN}, \c{S}iar and Keskin solved the same type equation, instead of taking sum, they considered the difference of two Fibonacci numbers and found solutions to the equation $F_{n}-F_{m}=2^{a}.$
Chim and Ziegler \cite{CHIM} considered the equations $F_{n_{1}}+F_{n_{2}}=2^{a_{1}}+2^{a_{2}}+2^{a_{3}}$ and $F_{m_{1}}+F_{m_{2}}+F_{m_{3}}=2^{t_{1}}+2^{t_{2}}$ and proved that $\max \{n_{1},n_{2},a_{1},a_{2},a_{3}\}\leq 18$ and $\max \{m_{1},m_{2},m_{3},t_{1},t_{2}\}\leq 16$, respectively.

The authors used lower bounds for linear forms in logarithms and a version of Baker-Davenport reduction method as their main tools to solve all the problems stated above. A natural question arises: What will be the solution if we replace Fibonacci numbers by balancing numbers? Therefore, in this note, we look at the Diophantine equation
 \begin{align}\label{eq:1}
 B_{n_{1}}+B_{n_{2}}=2^{a_{1}}+2^{a_{2}}+2^{a_{3}},
 \end{align}
where $B_{n}$ is the $n$-th balancing number with $n_{1}\geq n_{2}\geq  0$ and $a_{1}\geq a_{2} \geq a_{3} \geq 0$ and try to find all solutions using the same techniques.

The main result of this article is the following.
\begin{theorem}\label{th:1}
All non-negative integer solutions $(n_{1},n_{2},a_{1},a_{2},a_{3})$ of the equation \eqref{eq:1} are given by
\begin{multline*}
(n_{1},n_{2},a_{1},a_{2},a_{3})\in
\big\{(2,0,1,1,1), (2,0,2,0,0), (2,1,2,1,0),
(2,2,2,2,2), (2,2,3,1,1),\\
 (3,0,5,1,0), (3,1,4,4,2), (3,1,5,1,1), (3,2,5,3,0), (3,3,6,2,1)\big\}.
\end{multline*}
\end{theorem}
For the proof of Theorem \ref{th:1}, we run a program in \textit{Mathematica} and search all solutions $(n_{1},n_{2},a_{1},a_{2},a_{3})$ with $n_{1}<100$ to the equation \eqref{eq:1}. Then, we take $n_{1}>100$ and write \eqref{eq:1} in six different ways. We apply lower bounds for linear forms in logarithms to obtain an upper bound on $n_{1}=\max\{n_{1},n_{2},a_{1},a_{2},a_{3}\}.$ This is done in the following seven steps:

$Step~1$:  We find an upper bound
\begin{align*}
\min\left\{(a_{1}-a_{2})\log2, (n_{1}-n_{2})\log\alpha\right\}<8.22 \cdot 10^{12}(1+\log n_{1}).
\end{align*}
So, we divide into two cases:

~~~~~~~Case 1: $\min\{(a_{1}-a_{2})\log2, (n_{1}-n_{2})\log\alpha\}=(a_{1}-a_{2})\log2$

~~~~~~~Case 2: $\min\{(a_{1}-a_{2})\log2, (n_{1}-n_{2})\log\alpha\}=(n_{1}-n_{2})\log\alpha.$

$Step~2$: We consider case $1$ and show that
\begin{align*}
   \min\{(a_{1}-a_{3})\log2,(n_{1}-n_{2})\log\alpha\}<4 \cdot 10^{25}(1+\log n_{1})^{2}.
\end{align*}
We further divide case $1$ into two following sub-cases:

~~~~~~~Case 1A: $ \min\{(a_{1}-a_{3})\log2,(n_{1}-n_{2})\log\alpha\}=(a_{1}-a_{3})\log2$

~~~~~~~Case 1B: $ \min\{(a_{1}-a_{3})\log2,(n_{1}-n_{2})\log\alpha\}=(n_{1}-n_{2})\log\alpha.$

$Step ~3$: We consider case 1A and show that
\begin{align*}
   (n_{1}-n_{2})\log\alpha<2 \cdot 10^{38}(1+\log n_{1})^{3}.
\end{align*}

$Step~ 4$: We consider case 1B and show that
\begin{align*}
   (a_{1}-a_{3})\log2<9.96 \cdot 10^{37}(1+\log n_{1})^{3}.
\end{align*}

$Step~ 5$: We consider case 2 and show that
\begin{align*}
   (a_{1}-a_{2})\log2<2 \cdot 10^{25}(1+\log n_{1})^{2}.
\end{align*}

$Step~ 6$: We continue to consider case 2 and show that
\begin{align*}
   (a_{1}-a_{3})\log2<9.96 \cdot 10^{37}(1+\log n_{1})^{3}.
\end{align*}

$Step ~7$: Using the upper bounds $(a_{1}-a_{2})\log2, (a_{1}-a_{3})\log2, (n_{1}-n_{2})\log\alpha,$ we obtain an absolute upper bound for $n_{1}$ as
\begin{align*}
n_{1}<7.9 \cdot 10^{59}.
\end{align*}

We repeat all seven steps after finding an upper bound for $n_{1},$  but instead of lower bounds for linear forms in logarithms, we apply the Baker-Davenport reduction method. As a result, we have small absolute bounds and get to $n_{1}<86,$ a contradiction. In this way, we complete the proof of our main result.

In order to prove Theorem \ref{th:1}, we need some preliminary results which are discussed in the next section.
\section{Preliminaries}
Baker's theory of linear forms in logarithms of algebraic numbers plays an important role while solving various Diophantine equations. Here, we use several times the same to solve the equation \eqref{eq:1}, but before that, we recall some basic notations and results from algebraic number theory.

Let $\eta$ be an algebraic number with minimal primitive polynomial
\begin{align*}
     f\left(X\right) = a_{0}(X-\eta^{(1)})\ldots(X-\eta^{(k)}) \in {\mathbb{Z}\left[X\right]},
\end{align*}
where $a_{0}>0,$ and $\eta^{(i)}$'s are conjugates of $\eta.$ Then, the \textit{logarithmic height} of $\eta$ is defined by
\begin{align*}
    h(\eta)=\frac{1}{k}\left(\log a_{0}+\sum_{j=1}^{k}\max\{0,\log\vert\eta^{(j)}\vert\}\right).
\end{align*}
If $\eta=a/b$ is a rational number with gcd$(a,b)=1$ and $b>1,$ then $h(\eta)=\log (\max\{|a|,b\}).$ The following are some known properties of logarithmic height function:\begin{align*}
&h(\eta+\gamma)\leq h(\eta)+h(\gamma)+\log2,\\
&h(\eta\gamma^{\pm1})\leq h(\eta)+h(\gamma),\\
&h(\eta^{k})=|k|h(\eta),~~k\in\mathbb{Z}.
\end{align*}
The following theorem is a modified version of a result of Matveev (see \cite{Matveev} or \cite[Theorem~9.4]{Bugeaud1}) which provides a large upper bound for $n_{1}$ in \eqref{eq:1}.
 \begin{theorem}\label{th:2}
Let $\mathbb{L}$ be an algebraic number field of degree $d_{\mathbb{L}}.$ Let $\eta_{1},\eta_{2},\ldots , \eta_{l}\in \mathbb{L}$ be positive real numbers and $b_{1},b_{2},\ldots, b_{l}$ be nonzero integers. If
 $\Gamma=\prod_{i=1}^{l}\eta_{i}^{b_{i}}-1
$ is not zero, then
\begin{align*}
    \log \vert\Gamma\vert>-1.4 \cdot 30^{l+3} \cdot  l^{4.5} \cdot  d_{\mathbb{L}}^{2}(1+\log d_{\mathbb{L}})(1+\log D)A_{1}A_{2}\ldots A_{l},
\end{align*} where $D\geq\max\{\vert b_{1} \vert,\vert b_{2} \vert,\ldots, \vert b_{l} \vert\}$ and
$A_{1},A_{2},\ldots, A_{l}$ are positive real numbers such that \begin{align*}
    A_{j}\geq\max\{d_{\mathbb{L}}h\left( \eta_{j}\right),\vert\log\eta_{j}\vert, 0.16\}  ~~ \text{for}~   j=1,\ldots,l.
\end{align*}
\end{theorem}

We use the following method of Baker-Davenport due to Dujella and Peth\H o \cite{Dujella} to reduce the bound on $n_{1}.$
\begin{lemma}[\cite{Dujella}] \label{le:2.2}
Let M be a positive integer and $p/q$ be a convergent of the continued fraction of the irrational number $\tau$ such that $q>6M$. Let A, B, $\mu$ be some real numbers with $A>0$ and $B> 1.$ Let $\varepsilon:=\|\mu q\|-M\|\tau q \|, $ where $\|.\|$ denotes the distance from the nearest integer. If $\varepsilon>0,$ then there exists no solution to the inequality
\begin{align*}
    0<\lvert u\tau-v+\mu\rvert<AB^{-w},
\end{align*}
in positive integers $u,v, w$ with
\begin{align*}
    u\leq M~ and~ w\geq\frac{\log(A q/\varepsilon)}{\log B}.
\end{align*}
\end{lemma}
The following results will also be used to prove Theorem \ref{th:1}.
\begin{lemma}[\cite{GUZMAN}] \label{le:q}
Let $r\geq 1$ and $H>0$ be such that $H>(4r^{2})^r$ and $H>L/(\log L)^r.$ Then
\begin{align*}
L<2^{r}H(\log H)^r.
\end{align*}
\end{lemma}
\begin{lemma}\label{le:1}
All solutions of \eqref{eq:1} satisfy
 $(n_{1}-1)<\frac{\log3}{\log\alpha}+a_{1}\frac{\log2}{\log\alpha}$ and $n_{1}>(a_{1}-1)\frac{\log2}{\log\alpha}.$
\begin{proof}
From  \eqref{eq:1**} and \eqref{eq:1} we have
\begin{align*}
\alpha^{n_{1}-1}<B_{n_{1}}\leq B_{n_{1}}+B_{n_{2}}=2^{a_{1}}+2^{a_{2}}+2^{a_{3}}\leq 3 \cdot  2^{a_{1}}.
\end{align*}
Taking logarithm on both sides, we get
\[
\left(n_{1}-1\right)\log\alpha<\log 3 + a_{1} \log 2,
\]
which implies
\begin{align*}
    \left(n_{1}-1\right) < \frac{\log3}{\log\alpha}+a_{1}\frac{\log2}{\log\alpha}.
\end{align*}
On the other hand,
$2\alpha^{n_{1}}>2B_{n_{1}}\geq B_{n_{1}}+B_{n_{2}}=2^{a_{1}}+2^{a_{2}}+2^{a_{3}}>2^{a_{1}}$.
Taking logarithm on both sides, we get
\begin{align*}
    \log 2+ n_{1}\log\alpha>a_{1}\log2,
\end{align*} which implies
\[
n_{1}>(a_{1}-1)\frac{\log2}{\log\alpha}.
\qedhere
\]
\end{proof}
\end{lemma}
\section{Proof of Theorem \ref{th:1}}
Consider the Diophantine equation
\begin{align*}
B_{n_{1}}+B_{n_{2}}=2^{a_{1}}+2^{a_{2}}+2^{a_{3}}.
\end{align*}
First, we search the solutions to the above equation using \textit{Mathematica} for $n_{1}\leq 100$. Using Lemma \ref{le:1}, we calculate $a_{1}\leq256.$ By \textit{Mathematica}, for $0\leq n_{2}\leq n_{1}\leq 100$ and $0 \leq a_{3}\leq a_{2} \leq a_{1} \leq 256$, we find all the solutions that are listed in Theorem \ref{th:1}. Now, assume that $n_{1}>100.$
\subsection{An upper bound on $n_{1}$}
Using Binet's formula \eqref{eq:1} can be written as
\begin{equation}\label{eq:2}
   \frac{\alpha^{n_{1}}-\beta^{n_{1}}}{4\sqrt{2}}+\frac{\alpha^{n_{2}}-\beta^{n_{2}}}{4\sqrt{2}}=2^{a_{1}}+2^{a_{2}}+2^{a_{3}}.
\end{equation}
We write \eqref{eq:2} in the following six different ways and examine each one to prove our result.
\begin{equation}\label{eq:2a}
    \frac{\alpha^{n_{1}}}{4\sqrt{2}}-2^{a_{1}}=2^{a_{2}}+2^{a_{3}}+\frac{\beta^{n_{1}}}{4\sqrt{2}}-\frac{\alpha^{n_{2}}-\beta^{n_{2}}}{4\sqrt{2}}.
\end{equation}
\begin{equation}\label{eq:2b}
       \frac{\alpha^{n_{1}}}{4\sqrt{2}}-2^{a_{1}}-2^{a_{2}}=2^{a_{3}}+\frac{\beta^{n_{1}}}{4\sqrt{2}}-\frac{\alpha^{n_{2}}-\beta^{n_{2}}}{4\sqrt{2}}
\end{equation}
\begin{equation}\label{eq:2c}
       \frac{\alpha^{n_{1}}}{4\sqrt{2}}-2^{a_{1}}-2^{a_{2}}-2^{a_{3}}=\frac{\beta^{n_{1}}}{4\sqrt{2}}-\frac{\alpha^{n_{2}}-\beta^{n_{2}}}{4\sqrt{2}}
\end{equation}
\begin{equation}\label{eq:2d}
       \frac{\alpha^{n_{1}}}{4\sqrt{2}}+\frac{\alpha^{n_{2}}}{4\sqrt{2}}-2^{a_{1}}-2^{a_{2}}=2^{a_{3}}+\frac{\beta^{n_{1}}}{4\sqrt{2}}+\frac{\beta^{n_{2}}}{4\sqrt{2}}
\end{equation}
\begin{equation}\label{eq:2e}
       \frac{\alpha^{n_{1}}}{4\sqrt{2}}+\frac{\alpha^{n_{2}}}{4\sqrt{2}}-2^{a_{1}}=2^{a_{2}}+2^{a_{3}}+\frac{\beta^{n_{1}}}{4\sqrt{2}}+\frac{\beta^{n_{2}}}{4\sqrt{2}}
\end{equation}
\begin{equation}\label{eq:2f}
       \frac{\alpha^{n_{1}}}{4\sqrt{2}}+\frac{\alpha^{n_{2}}}{4\sqrt{2}}-2^{a_{1}}-2^{a_{2}}-2^{a_{3}}=\frac{\beta^{n_{1}}}{4\sqrt{2}}+\frac{\beta^{n_{2}}}{4\sqrt{2}}
\end{equation}

$\textbf{Step~1}$: First, we consider (\ref{eq:2a}). Here, we assume $n_{1}$ and $a_{1}$ to be large and collect the large terms involving $n_{1}$ and $a_{1}$ on the left side.
Taking absolute values on both sides of (\ref{eq:2a}), we get
\begin{align*}
\left\lvert\frac{\alpha^{n_{1}}}{4\sqrt{2}}-2^{a_{1}}\right\rvert&<2^{a_{2}+1}+\frac{\alpha^{n_{2}}}{4\sqrt{2}}+0.1\\
&<2.5\max\left\{2^{a_{2}},\alpha^{n_{2}}\right\}.
\end{align*}
Dividing both sides by $2^{a_{1}},$ we get
\begin{align*}
   \left\lvert\frac{\alpha^{n_{1}}}{4\sqrt{2}}2^{-a_{1}}-1\right\rvert<\max\left\{2.5 \cdot 2^{a_{2}-a_{1}},\frac{2.5\alpha^{n_{2}}}{2^{a_{1}}}\right\}<\max\left\{2.5 \cdot 2^{a_{2}-a_{1}},\frac{7.5\alpha^{n_{2}}}{\alpha^{n_{1}-1}}\right\}.
\end{align*}
Hence, we obtain
\begin{align}\label{eq:3}
   \left\lvert\frac{\alpha^{n_{1}}}{4\sqrt{2}}2^{-a_{1}}-1\right\rvert<43.72\max\left\{2^{a_{2}-a_{1}},\alpha^{n_{2}-n_{1}}\right\}.
\end{align}
Put \begin{align}\label{eq:4}
    \Gamma=\frac{\alpha^{n_{1}}}{4\sqrt{2}}2^{-a_{1}}-1.
\end{align}
Suppose $\Gamma=0,$ then $\alpha^{2n_{1}}\in \mathbb{Q}$ which is not possible for any $n_{1}>0.$ Therefore, $\Gamma\neq 0.$
To apply Theorem \ref{th:2} in \eqref{eq:4}, let
\begin{align*}
    \eta_{1}=\alpha,~ \eta_{2}=2,~  \eta_{3}=4\sqrt{2},~ b_{1}=n_{1},~ b_{2}=-a_{1},~ b_{3}=-1,~ l=3,
\end{align*}
where $\eta_{1},~\eta_{2},~\eta_{3}~\in{\mathbb{Q}(\alpha)}$ and $ b_{1},~ b_{2},~ b_{3}\in\mathbb{Z}.$
The degree $d_{\mathbb{L}}=[\mathbb{Q}(\alpha):\mathbb Q]$ is $2.$

Since $n_{1}>a_{1}>1,$ therefore $D =\max\{1,n_{1},\lvert a_{2}\rvert \}= n_{1}.$ We calculate the logarithmic heights of $\eta_{1}, \eta_{2}, \eta_{3}$ as follows:
\begin{align*}
    h(\eta_{1})=h\left(\alpha\right)=\frac{\log\alpha}{2},~h(\eta_{2})=\log2~ \text{and}~ h(\eta_{3})=\log(4\sqrt{2}).
\end{align*}
Thus, we can take
\begin{align*}
    A_{1}=\log\alpha,~~ A_{2}=2\log2~~ \text{and}~~ A_{3}=2\log(4\sqrt{2}).
\end{align*}
Applying Theorem \ref{th:2} we find
\begin{align*}
      \log \vert\Gamma\vert>-1.4 \cdot 30^{6} \cdot 3^{4.5} \cdot 2^2(1+\log2)(1+\log n_{1})(\log\alpha)(2\log2)(2\log(4\sqrt{2})).
\end{align*}
Comparing the above inequality with \eqref{eq:3} gives
\begin{align*}
\min\left\{(a_{1}-a_{2})\log2, (n_{1}-n_{2})\log\alpha\right\}<8.22 \cdot 10^{12}(1+\log n_{1}).
\end{align*}
Now, we divide into two cases.

\textbf{Case 1:} $\min\{(a_{1}-a_{2})\log2, (n_{1}-n_{2})\log\alpha\}=(a_{1}-a_{2})\log2$.

\textbf{Case 2:} $\min\{(a_{1}-a_{2})\log2, (n_{1}-n_{2})\log\alpha\}=(n_{1}-n_{2})\log\alpha$.

$\textbf{Step~2}$: First, we consider case 1 and assume that
\begin{align}\label{in:4}
\min\{(a_{1}-a_{2})\log2, (n_{1}-n_{2})\log\alpha\}=(a_{1}-a_{2})\log2<8.22 \cdot 10^{12}(1+\log n_{1}).
\end{align}
Assuming $n_{1}, a_{1}$ and $a_{2}$ to be large and collecting large terms on the left hand side, we consider \eqref{eq:2b}. Taking absolute values on both sides of \eqref{eq:2b}, we have
\begin{equation*}
       \left\lvert\frac{\alpha^{n_{1}}}{4\sqrt{2}}-2^{a_{1}}-2^{a_{2}}\right\rvert=\left\lvert2^{a_{3}}+\frac{\beta^{n_{1}}}{4\sqrt{2}}-\frac{\alpha^{n_{2}}-\beta^{n_{2}}}{4\sqrt{2}}\right\rvert,
\end{equation*}
which implies
\begin{align*}
\left\lvert \frac{\alpha^{n_{1}}}{4\sqrt{2}}-2^{a_{1}}-2^{a_{2}}\right\rvert<2^{a_{3}}+\frac{\alpha^{n_{2}}}{4\sqrt{2}}+0.1<1.2\max\{2^{a_{3}},\alpha^{n_{2}}\}.
\end{align*}
Dividing both sides by $\frac{\alpha^{n_{1}}}{4\sqrt{2}},$ we obtain
\begin{align*}
     \left\lvert 1-\alpha^{-n_{1}}2^{a_{2}}4\sqrt{2}(2^{a_{1}-a_{2}}+1)\right\rvert&<\max\left\{\frac{(1.2)(4\sqrt{2})}{\alpha^{n_{1}}} \cdot 2^{a_{3}},(1.2)(4\sqrt{2})\alpha^{n_{2}-n_{1}}\right\}\\
     &\leq \max\left\{\frac{(1.2)(4\sqrt{2})}{2^{a_{1}-1}} \cdot 2^{a_{3}},(1.2)(4\sqrt{2})\alpha^{n_{2}-n_{1}}\right\}.
\end{align*}
Hence, we obtain
\begin{align}\label{in:3}
\left\lvert 1-\alpha^{-n_{1}}2^{a_{2}}4\sqrt{2}(2^{a_{1}-a_{2}}+1)\right\rvert<13.57\max\left\{{2^{a_{3}-a_{1}}},\alpha^{n_{2}-n_{1}}\right\}.
\end{align}
Put
\begin{align*}
    \Gamma_{1}=1-\alpha^{-n_{1}}2^{a_{2}}4\sqrt{2}(2^{a_{1}-a_{2}}+1).
\end{align*}
By similar arguments as before we can show that $\Gamma_{1}\neq0.$
With the notations of Theorem~\ref{th:2}, we take
\begin{align*}
    \eta_{1}=\alpha,~ \eta_{2}=2,~  \eta_{3}=4\sqrt{2}(2^{a_{1}-a_{2}}+1),~ b_{1}=-n_{1},~ b_{2}=a_{2},~ b_{3}=1,~ l=3.
\end{align*}
Since $a_{2}<n_{1},$ we take $D=n_{1}.$ As before, we have the same logarithmic heights for $\eta_{1}$ and  $\eta_{2}.$ Thus $A_{1}$ and $A_{2}$ remain unchanged. Computing the height of $\eta_{3},$ we have
   \begin{align*}
   h(\eta_{3})&= h(4\sqrt{2}(2^{a_{1}-a_{2}}+1))\\
  &\leq h(4\sqrt{2})+h(2^{a_{1}-a_{2}}+1)\\
  & \leq \log(4\sqrt{2})+(a_{1}-a_{2})\log2+\log2.
   \end{align*}
Hence, from \eqref{in:4}, we get
\begin{align*}
 h(\eta_{3})<8.23 \cdot 10^{12}(1+\log n_{1}).
\end{align*}
So, we take
\begin{align*}
    A_{3}=16.46 \cdot 10^{12}(1+\log n_{1}).
\end{align*}
Using all these values in Theorem \ref{th:2}, we have
\begin{align*}
      \log \vert\Gamma_{1}\vert > -1.4 \cdot 30^{6} \cdot 3^{4.5} \cdot 2^2(1+\log2)(1+\log n_{1})(\log\alpha)(2\log2)(16.46 \cdot 10^{12}(1+\log n_{1})).
\end{align*}
Comparing the above inequality with \eqref{in:3} gives
\begin{align*}
   \min\{(a_{1}-a_{3})\log2,(n_{1}-n_{2})\log\alpha\}<4 \cdot 10^{25}(1+\log n_{1})^{2}.
\end{align*}
Now, we divide this into two sub-cases.

\textbf{Case 1A:} $ \min\{(a_{1}-a_{3})\log2,(n_{1}-n_{2})\log\alpha\}=(a_{1}-a_{3})\log2$.

\textbf{Case 1B:} $ \min\{(a_{1}-a_{3})\log2,(n_{1}-n_{2})\log\alpha\}=(n_{1}-n_{2})\log\alpha$.

$\textbf{Step~3}$: Assume the first sub-case, that is
\begin{align}\label{en:Z}
\min\{(a_{1}-a_{3})\log2,(n_{1}-n_{2})\log\alpha\}=(a_{1}-a_{3})\log2<4 \cdot 10^{25}(1+\log n_{1})^{2}.
\end{align}
In this step, we consider $n_{1}, a_{1}, a_{2}$ and $a_{3}$ to be large. By collecting large terms on the left side, we consider (\ref{eq:2c}), that is
\begin{equation*}
       \left\lvert\frac{\alpha^{n_{1}}}{4\sqrt{2}}-2^{a_{1}}-2^{a_{2}}-2^{a_{3}}\right\rvert=\left\lvert\frac{\beta^{n_{1}}}{4\sqrt{2}}-\frac{\alpha^{n_{2}}-\beta^{n_{2}}}{4\sqrt{2}}\right\rvert,
\end{equation*}
which implies
\begin{align*}
\left\lvert \frac{\alpha^{n_{1}}}{4\sqrt{2}}-2^{a_{1}}-2^{a_{2}}-2^{a_{3}}\right\rvert<\frac{\alpha^{n_{2}}}{4\sqrt{2}}+0.1<0.3\alpha^{n_{2}}.
\end{align*}
Dividing both sides by $\frac{\alpha^{n_{1}}}{4\sqrt{2}}$, we obtain
\begin{align}\label{en:B}
     \left\lvert 1-\alpha^{-n_{1}}2^{a_{1}}4\sqrt{2}(1+2^{a_{2}-a_{1}}+2^{a_{3}-a_{1}})\right\rvert&<0.3\alpha^{n_{2}}\left(\frac{4\sqrt{2}}{\alpha^{n_{1}}}\right)=1.7\alpha^{n_{2}-n_{1}}.
\end{align}
Put
\begin{align*}
    \Gamma_{A}=1-\alpha^{-n_{1}}2^{a_{1}}4\sqrt{2}(1+2^{a_{2}-a_{1}}+2^{a_{3}-a_{1}}).
\end{align*}
We can show that $\Gamma_{A}\neq0$.
Take
\begin{align*}
    \eta_{1}=\alpha,~ \eta_{2}=2,~  \eta_{3}=4\sqrt{2}(1+2^{a_{2}-a_{1}}+2^{a_{3}-a_{1}}),~ b_{1}=-n_{1},~ b_{2}=a_{1},~ b_{3}=1.
\end{align*}
Computing the logarithmic height of $\eta_{3},$ we get
   \begin{align*}
   h(\eta_{3})&= h(4\sqrt{2}(1+2^{a_{2}-a_{1}}+2^{a_{3}-a_{1}}))\\
  &\leq h(4\sqrt{2})+h(1+2^{a_{2}-a_{1}}+2^{a_{3}-a_{1}})\\
  & \leq \log(4\sqrt{2})+(a_{1}-a_{2})\log2+(a_{1}-a_{3})\log2+2\log2.
   \end{align*}
Hence, from \eqref{in:4} and \eqref{en:Z}, we get
\begin{align*}
 h(\eta_{3})<4.1 \cdot 10^{25}(1+\log n_{1})^{2}.
\end{align*}
So, we take
\begin{align*}
     A_{3}=8.2 \cdot 10^{25}(1+\log n_{1})^{2}.
\end{align*}
The parameters $A_{1}$ and $A_{2}$ remain unchanged as before.
Using all these values in Theorem~\ref{th:2}, we have
\begin{align*}
      \log \vert\Gamma_{A}\vert>-1.4 \cdot 30^{6} \cdot 3^{4.5} \cdot 2^2(1+\log2)(1+\log n_{1})(\log\alpha)(2\log2)(8.2 \cdot 10^{25}(1+\log n_{1})^{2}).
\end{align*}
Comparing the above inequality with \eqref{en:B} gives
\begin{align*}
   (n_{1}-n_{2})\log\alpha<2 \cdot 10^{38}(1+\log n_{1})^{3}.
\end{align*}

$\textbf{Step~4}$: Now, we consider the second sub-case, that is
 \begin{align}\label{en:A1}
\min\{(a_{1}-a_{3})\log2,(n_{1}-n_{2})\log\alpha\}=(n_{1}-n_{2})\log\alpha<4 \cdot 10^{25}(1+\log n_{1})^{2}.
\end{align}
Equation \eqref{eq:2d} implies
\begin{align*}
\left\lvert \frac{\alpha^{n_{2}}(1+\alpha^{n_{1}-n_{2}})}{4\sqrt{2}}-2^{a_{2}}(2^{a_{1}-a_{2}}+1)\right\rvert<1.1 \cdot 2^{a_{3}}.
\end{align*}
Dividing both sides by $2^{a_{2}}(2^{a_{1}-a_{2}}+1)$, we obtain
\begin{align}\label{en:B1}
     \left\lvert \alpha^{n_{2}}2^{-a_{2}}\frac{(1+\alpha^{n_{1}-n_{2}})}{4\sqrt{2}(2^{a_{1}-a_{2}}+1)}-1\right\rvert&<1.1 \cdot 2^{a_{3}-a_{1}}.
\end{align}
Take
\begin{align*}
    \Gamma_{B}= \alpha^{n_{2}}2^{-a_{2}}\frac{(1+\alpha^{n_{1}-n_{2}})}{4\sqrt{2}(2^{a_{1}-a_{2}}+1)}-1,
\end{align*}
with $\eta_{1}=\alpha,~ \eta_{2}=2,~  \eta_{3}=\frac{(1+\alpha^{n_{1}-n_{2}})}{4\sqrt{2}(2^{a_{1}-a_{2}}+1)},~ b_{1}=n_{2},~ b_{2}=-a_{2},~ b_{3}=1$.
Since $a_{2}<n_{2}<n_{1},$ $D=n_{1}.$ The height of $\eta_{3}$ is calculated as
\begin{align*}
   h(\eta_{3})&= h\left(\frac{(1+\alpha^{n_{1}-n_{2}})}{4\sqrt{2}(2^{a_{1}-a_{2}}+1)}\right)\\
  &\leq h(1+\alpha^{n_{1}-n_{2}})+h(4\sqrt{2}(2^{a_{1}-a_{2}}+1))\\
  & \leq (n_{1}-n_{2})h(\alpha)+h(4\sqrt{2})+(a_{1}-a_{2})h(2)+2\log2\\
  &=(n_{1}-n_{2})\frac{\log\alpha}{2}+\log(4\sqrt{2})+(a_{1}-a_{2})\log2+2\log2.
 \end{align*}
Hence, from \eqref{in:4} and \eqref{en:A1}, we get
\begin{align*}
 h(\eta_{3})<2.1 \cdot 10^{25}(1+\log n_{1})^{2}.
\end{align*}
So, we take
\begin{align*}
    A_{3}=4.2 \cdot 10^{25}(1+\log n_{1})^{2}.
\end{align*}
Applying Theorem \ref{th:2}, we have
\begin{align*}
      \log \vert\Gamma_{B}\vert>-1.4 \cdot 30^{6} \cdot 3^{4.5} \cdot 2^2(1+\log2)(1+\log n_{1})(\log\alpha)(2\log2)(4.2 \cdot 10^{25}(1+\log n_{1})^{2}).
\end{align*}
Comparing the above inequality with \eqref{en:B1} gives
\begin{align*}
   (a_{1}-a_{3})\log2<9.96 \cdot 10^{37}(1+\log n_{1})^{3}.
\end{align*}

$\textbf{Step~5}$: Now, we consider case $2,$ that is
\begin{align}\label{in:MM}
\min\{(a_{1}-a_{2})\log2, (n_{1}-n_{2})\log\alpha\}=(n_{1}-n_{2})\log\alpha<8.22 \cdot 10^{12}(1+\log n_{1}).
\end{align}
Equation \eqref{eq:2e} implies
\begin{align*}
\left\lvert \frac{\alpha^{n_{2}}(1+\alpha^{n_{1}-n_{2}})}{4\sqrt{2}}-2^{a_{1}}\right\rvert<2.2 \cdot 2^{a_{2}}.
\end{align*}
Dividing both sides by $2^{a_{1}},$ we obtain
\begin{align}\label{en:B2}
     \left\lvert \alpha^{n_{2}}2^{-a_{1}}\frac{(1+\alpha^{n_{1}-n_{2}})}{4\sqrt{2}}-1\right\rvert&<2.2 \cdot 2^{a_{2}-a_{1}}.
\end{align}
Put
\begin{align*}
     \Gamma_{2}= \alpha^{n_{2}}2^{-a_{1}}\frac{(1+\alpha^{n_{1}-n_{2}})}{4\sqrt{2}}-1.
\end{align*}
We can show that $\Gamma_{2}\neq0.$
 With the notations of Theorem \ref{th:2}, we take
\begin{align*}
    \eta_{1}=\alpha,~ \eta_{2}=2,~  \eta_{3}=\frac{(1+\alpha^{n_{1}-n_{2}})}{4\sqrt{2}},~ b_{1}=n_{2},~ b_{2}=-a_{1},~ b_{3}=1.
\end{align*}
Since $a_{2}<n_{2}<n_{1},$ $D=n_{1}.$ Computing the logarithmic height of $\eta_{3},$ we get
   \begin{align*}
   h(\eta_{3})&= h\left(\frac{1+\alpha^{n_{1}-n_{2}}}{4\sqrt{2}}\right)\\
  &\leq h(1+\alpha^{n_{1}-n_{2}})+h(4\sqrt{2})\\
  & \leq (n_{1}-n_{2})h(\alpha)+h(4\sqrt{2})+\log2\\
  &=(n_{1}-n_{2})\frac{\log\alpha}{2}+\log(4\sqrt{2})+\log2.
   \end{align*}
Hence, from \eqref{in:MM}, we obtain
\begin{align*}
 h(\eta_{3})<4.12 \cdot 10^{12}(1+\log n_{1}).
\end{align*}
So, we take
\begin{align*}
    A_{3}=8.24 \cdot 10^{12}(1+\log n_{1}).
\end{align*} The value of $A_{1}$ and $A_{2}$ remain same as before.
Applying Theorem \ref{th:2}, we have
\begin{align*}
      \log \vert\Gamma_{2}\vert>-1.4 \cdot 30^{6} \cdot 3^{4.5} \cdot 2^2(1+\log2)(1+\log n_{1})(\log\alpha)(2\log2)(8.24 \cdot 10^{12}(1+\log n_{1})).
\end{align*}
Comparing the above inequality with \eqref{en:B2} gives
\begin{align}\label{in:a3}
   (a_{1}-a_{2})\log2<2 \cdot 10^{25}(1+\log n_{1})^{2}.
\end{align}

$\textbf{Step~6}$: We apply Theorem \ref{th:2} once more to obtain an upper bound for $(a_{1}-a_{3})\log2.$ The derivation is similar to case 1B.
By the similar derivation as in step $4$, we obtain
\begin{align}\label{en:B12}
     \left\lvert \alpha^{n_{2}}2^{-a_{2}}\frac{(1+\alpha^{n_{1}-n_{2}})}{4\sqrt{2}(2^{a_{1}-a_{2}}+1)}-1\right\rvert&<1.1 \cdot 2^{a_{3}-a_{1}}.
\end{align}
We estimate the height of $\eta_{3}$ as
\begin{align*}
   h(\eta_{3})&= h\left(\frac{(1+\alpha^{n_{1}-n_{2}})}{4\sqrt{2}(2^{a_{1}-a_{2}}+1)}\right)\\
  & \leq (n_{1}-n_{2})h(\alpha)+h(4\sqrt{2})+(a_{1}-a_{2})h(2)+2\log2\\
  &=(n_{1}-n_{2})\frac{\log\alpha}{2}+\log(4\sqrt{2})+(a_{1}-a_{2})\log2+2\log2.
\end{align*}
Hence, from \eqref{in:MM} and \eqref{in:a3}, we get
\begin{align*}
 h(\eta_{3})<2.1 \cdot 10^{25}(1+\log n_{1})^{2}.
\end{align*}
So, we take
\begin{align*}
    A_{3}=4.2 \cdot 10^{25}(1+\log n_{1})^{2}.
\end{align*}
Applying Theorem \ref{th:2}, we have
\begin{align*}
      \log \vert\Gamma_{B}\vert>-1.4 \cdot 30^{6} \cdot 3^{4.5} \cdot 2^2(1+\log2)(1+\log n_{1})(\log\alpha)(2\log2)(4.2 \cdot 10^{25}(1+\log n_{1})^{2}).
\end{align*}
Comparing the above inequality with \eqref{en:B12} gives
\begin{align*}
   (a_{1}-a_{3})\log2<9.96 \cdot 10^{37}(1+\log n_{1})^{3}.
\end{align*}
We summarize our results obtained so far in the following table.
\smallskip\newline\medskip
\resizebox{\textwidth}{!}{%
\begin{tabular}{l|l|l|l}
Upper bound of & Case 1A & Case 1B & Case 2 \\
\hline
$(a_{1}-a_{2})\log2$ & $8.22 \cdot 10^{12}(1+\log n_{1})$ & $8.22 \cdot 10^{12}(1+\log n_{1})$ & $2 \cdot 10^{25}(1+\log n_{1})^{2}$\\
$(a_{1}-a_{3})\log2$ &
$4 \cdot 10^{25}(1+\log n_{1})^{2}$ & $9.96 \cdot 10^{37}(1+\log n_{1})^{3}$ & $9.96 \cdot 10^{37}(1+\log n_{1})^{3}$  \\
$(n_{1}-n_{2})\log\alpha$ & $2 \cdot 10^{38}(1+\log n_{1})^{3}$ &
$4 \cdot 10^{25}(1+\log n_{1})^{2}$ & $8.22 \cdot 10^{12}(1+\log n_{1})$
\end{tabular}
}

$\textbf{Step~7}$: Lastly, we consider \eqref{eq:2f}, that is
\begin{equation*}
       \frac{\alpha^{n_{1}}}{4\sqrt{2}}+\frac{\alpha^{n_{2}}}{4\sqrt{2}}-2^{a_{1}}-2^{a_{2}}-2^{a_{3}}=\frac{\beta^{n_{1}}}{4\sqrt{2}}+\frac{\beta^{n_{2}}}{4\sqrt{2}}.
\end{equation*}
Taking absolute values on both sides, we have
\begin{align*}
\left\lvert \frac{\alpha^{n_{1}}(1+\alpha^{n_{2}-n_{1}})}{4\sqrt{2}}-2^{a_{1}}(1+2^{a_{2}-a_{1}}+2^{a_{3}-a_{1}})\right\rvert<0.1.
\end{align*}
Dividing both sides by $\frac{\alpha^{n_{1}}(1+\alpha^{n_{2}-n_{1}})}{4\sqrt{2}}$ gives
\begin{align}\label{en:B21}
     \left\lvert 1-\alpha^{-n_{1}}2^{a_{1}}\frac{4\sqrt{2}(1+2^{a_{2}-a_{1}}+2^{a_{3}-a_{1}})}{(1+\alpha^{n_{2}-n_{1}})}\right\rvert&<0.6 \cdot \alpha^{-n_{1}}.
\end{align}
Put
\begin{align*}
    \Gamma_{3}=  \left\lvert 1-\alpha^{-n_{1}}2^{a_{1}}\frac{4\sqrt{2}(1+2^{a_{2}-a_{1}}+2^{a_{3}-a_{1}})}{(1+\alpha^{n_{2}-n_{1}})}\right\rvert.
\end{align*}
Using similar arguments as before we can show that $\Gamma_{3}\neq0.$
With the notations of Theorem~\ref{th:2}, we take
\begin{align*}
    \eta_{1}=\alpha,~ \eta_{2}=2,~  \eta_{3}=\frac{4\sqrt{2}(1+2^{a_{2}-a_{1}}+2^{a_{3}-a_{1}})}{(1+\alpha^{n_{2}-n_{1}})},~ b_{1}=-n_{1},~ b_{2}=a_{1},~ b_{3}=1.
\end{align*}
Since $a_{1}<n_{1},$ $D=n_{1}.$ Computing the logarithmic height of $\eta_{3},$ we get
   \begin{align*}
   h(\eta_{3})&= h\left(\frac{4\sqrt{2}(1+2^{a_{2}-a_{1}}+2^{a_{3}-a_{1}})}{(1+\alpha^{n_{2}-n_{1}})}\right)\\
  &\leq h(4\sqrt{2}(1+2^{a_{2}-a_{1}}+2^{a_{3}-a_{1}}))+h(1+\alpha^{n_{2}-n_{1}})\\
  & \leq h(4\sqrt{2})+(a_{1}-a_{2})h(2)+(a_{1}-a_{3})h(2)+(n_{1}-n_{2})h(\alpha)+3\log2\\
  &=\log(4\sqrt{2})+(a_{1}-a_{2})\log2+(a_{1}-a_{3})\log2+(n_{1}-n_{2})\frac{\log\alpha}{2}+3\log2.
   \end{align*}
Hence, we get
\begin{align*}
 h(\eta_{3})<9.97 \cdot 10^{37}(1+\log n_{1})^{3}.
\end{align*}
So, we take
\begin{align*}
  A_{3}=19.95 \cdot 10^{37}(1+\log n_{1})^{3}.
\end{align*}
Applying Theorem \ref{th:2}, we have
\[
    \log |\Gamma_{3}|
    >-1.4 \cdot 30^{6} \cdot 3^{4.5} \cdot 2^2(1+\log2)(1+\log n_{1})(\log\alpha)(2\log2)(19.95 \cdot 10^{37}(1+\log n_{1})^{3}).
\]
Comparing the above inequality with \eqref{en:B21} gives
\begin{align*}
   n_{1}\log\alpha<4.73 \cdot 10^{50}(1+\log n_{1})^{4}.
\end{align*}
With the notation of Lemma \ref{le:q}, we take $r=4$, $L=n$ and $H=\frac{4.73 \cdot 10^{50}}{\log\alpha}$. Applying the lemma, we have
\begin{align*}
n_{1}&<2^{4}\left(\frac{4.73 \cdot 10^{50}}{\log\alpha}\right)\left(\log\left(\frac{4.73 \cdot 10^{50}}{\log\alpha}\right)\right)^{4}\\
&<7.9 \cdot 10^{59}.
\end{align*}
The bound on $n_{1}$ is too large. So, in the next subsection, we reduce this bound using Lemma \ref{le:2.2}.
\subsection{Bound Reduction}
To reduce the bound on $n_{1},$ we use the following steps.

$\textbf{Step~1}$:
Put
\begin{align*}
    \Lambda=n_{1}\log\alpha-a_{1}\log 2-\log\left(4\sqrt{2}\right).
\end{align*}
The inequality \eqref{eq:3} can be written as
\begin{align*}
      \left\lvert\frac{\alpha^{n_{1}}}{4\sqrt{2}}2^{-a_{1}}-1\right\rvert=\lvert e^{\Lambda}-1\rvert<43.72\max\left\{2^{a_{2}-a_{1}},\alpha^{n_{2}-n_{1}}\right\}.
\end{align*}
Observe that $\Lambda\neq0$ as $e^{\Lambda}-1=\Gamma\neq0.$
Assuming $\min\left\{{a_{1}-a_{2}},{n_{1}-n_{2}}\right\}\geq7,$ the right-hand side in the above inequality is at most $\frac{1}{2}.$
The inequality $\lvert e^{z}-1\rvert<y$ for real values of $z$ and $y$ implies $z<2y.$ Thus, we get
\begin{align*}
     \lvert \Lambda\rvert<87.44\max\left\{2^{a_{2}-a_{1}},\alpha^{n_{2}-n_{1}}\right\},
\end{align*}
which implies that
\begin{align*}
    \left\lvert n_{1}\log\alpha-a_{1}\log 2-\log\left(4\sqrt{2}\right)\right\rvert<87.44\max\left\{2^{a_{2}-a_{1}},\alpha^{n_{2}-n_{1}}\right\}.
\end{align*}
Dividing both sides by $\log 2$ gives
\begin{align*}
\left\lvert n_{1}\left(\frac{\log\alpha}{\log 2}\right)-a_{1}+\frac{\log\left(1/4\sqrt{2}\right)}{\log 2}\right\rvert&<\max\left\{\frac{87.44}{\log2} \cdot 2^{a_{2}-a_{1}},\frac{87.44}{\log2}\alpha^{n_{2}-n_{1}}\right\}\\
&<\max\left\{127 \cdot 2^{-(a_{1}-a_{2})},127\alpha^{-(n_{1}-n_{2}})\right\}.
\end{align*}
We let
\begin{align*}
     &u=n_{1},~ \tau=\left(\frac{\log\alpha}{\log 2}\right),~ v=a_{1},~ \mu=\frac{\log\left(1/4\sqrt{2}\right)}{\log 2},
     ~\text{with}\\
     &\left(A,B,w\right)=(127,2,(a_{1}-a_{2}))~\text{or}~ (127,\alpha,(n_{1}-n_{2})).
\end{align*}
Choose $M=7.9 \cdot 10^{59}$. We find $q_{126}$ exceeds $6M$ with $\varepsilon=\|\mu q_{126}\|-M\|\tau q_{126} \|=0.5$. By virtue of Lemma \ref{le:2.2}, we get $a_{1}-a_{2}\leq 214$ or $n_{1}-n_{2}\leq 84$.
Now, we divide this into two cases.

\textbf{Case 1:} $a_{1}-a_{2}\leq 214$

\textbf{Case 2:} $n_{1}-n_{2}\leq 84$

$\textbf{Step~2}$: First, we consider case $1$.
Let
 \begin{align*}
    \Lambda_{1}=-n_{1}\log\alpha+a_{2}\log 2+\log\left(4\sqrt{2}(1+2^{a_{1}-a_{2}})\right).
\end{align*}
The inequality \eqref{in:3} can be written as
\begin{align*}
\left\lvert e^{\Lambda_{1}}-1\right\rvert=\lvert \Gamma_{1}\rvert<13.57\max\left\{{2^{a_{3}-a_{1}}},\alpha^{n_{2}-n_{1}}\right\}.
\end{align*}
Assuming $\min\left\{{a_{1}-a_{3}},{n_{1}-n_{2}}\right\}\geq5$, the right-hand side in the above inequality at most $\frac{1}{2}.$ Thus, we get
 \begin{align*}
    \left\lvert n_{1}\log\alpha-a_{2}\log 2+\log\left(1/(4\sqrt{2}(1+2^{a_{1}-a_{2}}))\right)\right\rvert<27.14\max\left\{{2^{a_{3}-a_{1}}},\alpha^{n_{2}-n_{1}}\right\}.
 \end{align*}
 Dividing both sides by $\log 2$ gives
 \begin{align*}
     \left\lvert n_{1}\left(\frac{\log\alpha}{\log 2}\right)-a_{2}+\frac{\log\left(1/(4\sqrt{2}(1+2^{a_{1}-a_{2}}))\right)}{\log 2}\right\rvert&<\max\left\{\frac{27.52}{\log2} \cdot 2^{a_{3}-a_{1}},\frac{27.52}{\log2}\alpha^{n_{2}-n_{1}}\right\}\\ &<\max\left\{40 \cdot 2^{-(a_{1}-a_{3})},40\alpha^{-(n_{1}-n_{2}})\right\}.
 \end{align*}
Let
 \begin{align*}
     u=n_{1},~ \tau=\left(\frac{\log\alpha}{\log 2}\right),~ v=a_{2},~ \mu=\frac{\log\left(1/(4\sqrt{2}(1+2^{a_{1}-a_{2}}))\right)}{\log 2},
 \end{align*}
with $\left(A,B,w\right)=(40,2,(a_{1}-a_{3}))~\text{or}~ (40,\alpha,(n_{1}-n_{2}))$.
With the same $M$, we find $q_{124}$ exceeds $6M$ with $\varepsilon>0.00179287.$ By virtue of Lemma \ref{le:2.2} for $(a_{1}-a_{2})\leq 214,$ we get $a_{1}-a_{3}\leq 218$ or $n_{1}-n_{2}\leq 86.$

Again, we divide case $1$ into two sub-cases.

\textbf{Case 1A:} $a_{1}-a_{3}\leq 218$

\textbf{Case 1B:} $n_{1}-n_{2}\leq 86$

$\textbf{Step~3}$: We consider 1A. Put
 \begin{align*}
    \Lambda_{A}=-n_{1}\log\alpha+a_{1}\log 2+\log\left(4\sqrt{2}(1+2^{a_{2}-a_{1}}+2^{a_{3}-a_{1}})\right).
\end{align*}
Then, inequality \eqref{en:B} can be written as
\begin{align*}
\left\lvert e^{\Lambda_{A}}-1\right\rvert=\lvert \Gamma_{A}\rvert<1.7\alpha^{n_{2}-n_{1}}.
\end{align*}
Assuming $(n_{1}-n_{2})\geq 1,$ we get
 \begin{align*}
    \left\lvert n_{1}\log\alpha-a_{1}\log 2+\log\left(1/(4\sqrt{2}(1+2^{a_{2}-a_{1}}+2^{a_{3}-a_{1}}))\right)\right\rvert<3.4\alpha^{n_{2}-n_{1}},
 \end{align*}
which implies
 \begin{align*}
     \left\lvert n_{1}\left(\frac{\log\alpha}{\log 2}\right)-a_{1}+\frac{\log\left(1/(4\sqrt{2}(1+2^{a_{2}-a_{1}}+2^{a_{3}-a_{1}}))\right)}{\log 2}\right\rvert&<\frac{3.4}{\log2}\alpha^{n_{2}-n_{1}}\\
     &<5\alpha^{-(n_{1}-n_{2})}.
 \end{align*}
 Let
 \begin{align*}
     u=n_{1},~ \tau=\left(\frac{\log\alpha}{\log 2}\right),~ v=a_{1},~ \mu=\frac{\log\left(1/(4\sqrt{2}(1+2^{a_{2}-a_{1}}+2^{a_{3}-a_{1}}))\right)}{\log 2},
 \end{align*}
with $\left(A,B,w\right)=(5,\alpha,(n_{1}-n_{2})).$
With the same $M$, we estimate $\varepsilon>0.0000354843$. Applying Lemma \ref{le:2.2} for $(a_{1}-a_{2})\leq 214$ and $(a_{1}-a_{3})\leq 218,$ we get $n_{1}-n_{2}\leq 87.$

$\textbf{Step~4}$: We consider the case 1B. Put
 \begin{align*}
    \Lambda_{B}=n_{2}\log\alpha-a_{2}\log 2+\log\frac{(1+\alpha^{n_{1}-n_{2}})}{4\sqrt{2}(2^{a_{1}-a_{2}}+1)}.
\end{align*}
The inequality \eqref{en:B1} can be written as
\begin{align*}
\left\lvert e^{\Lambda_{B}}-1\right\rvert=\lvert \Gamma_{B}\rvert<1.1 \cdot 2^{a_{3}-a_{1}}.
\end{align*}
Assuming $(a_{1}-a_{3})\geq 2$, we get
 \begin{align*}
    \left\lvert n_{2}\log\alpha-a_{2}\log 2+\log\frac{(1+\alpha^{n_{1}-n_{2}})}{4\sqrt{2}(2^{a_{1}-a_{2}}+1)}\right\rvert<2.2 \cdot 2^{-(a_{1}-a_{3})},
 \end{align*}
 which implies
 \begin{align*}
     \left\lvert n_{2}\left(\frac{\log\alpha}{\log 2}\right)-a_{2}+\frac{\log\left((1+\alpha^{n_{1}-n_{2}})/(4\sqrt{2}(2^{a_{1}-a_{2}}+1))\right)}{\log 2}\right\rvert&<\frac{2.2}{\log2} \cdot 2^{a_{3}-a_{1}}\\
     &<3.1 \cdot 2^{-({a_{1}-a_{3}})}.
 \end{align*}
Let
 \begin{align*}
     u=n_{2},~ \tau=\left(\frac{\log\alpha}{\log 2}\right),~ v=a_{2},~ \mu=\frac{\log\left((1+\alpha^{n_{1}-n_{2}})/(4\sqrt{2}(2^{a_{1}-a_{2}}+1))\right)}{\log 2},
 \end{align*}
with $\left(A,B,w\right)=(3.1,2,({a_{1}-a_{3}})).$
With the same $M$ we find $\varepsilon>0.0000119685$. Applying Lemma \ref{le:2.2} for $(a_{1}-a_{2})\leq 214$ and $(n_{1}-n_{2})\leq 86,$ we get $a_{1}-a_{3}\leq 222.$

$\textbf{Step~5}$: Now, consider case $2.$
Take
 \begin{align*}
    \Lambda_{2}=n_{2}\log\alpha-a_{1}\log 2+\log\frac{(1+\alpha^{n_{1}-n_{2}})}{4\sqrt{2}}.
\end{align*}
The inequality \eqref{en:B2} can be written as
\begin{align*}
\left\lvert e^{\Lambda_{2}}-1\right\rvert=\lvert \Gamma_{2}\rvert<2.2 \cdot 2^{a_{2}-a_{1}}.
\end{align*}
Assuming $(a_{1}-a_{2})\geq 3,$ we get
 \begin{align*}
    \left\lvert n_{2}\log\alpha-a_{1}\log 2+\log\frac{(1+\alpha^{n_{1}-n_{2}})}{4\sqrt{2}}\right\rvert<4.4 \cdot 2^{-(a_{1}-a_{2})}.
 \end{align*}
Dividing both sides by $\log 2$ gives
 \begin{align*}
     \left\lvert n_{2}\left(\frac{\log\alpha}{\log 2}\right)-a_{1}+\frac{\log\left((1+\alpha^{n_{1}-n_{2}})/4\sqrt{2}\right)}{\log 2}\right\rvert&<\frac{4.4}{\log2} \cdot 2^{-({a_{1}-a_{2}})}\\
     &<6.3 \cdot 2^{-({a_{1}-a_{2}})}.
 \end{align*}
 Let
 \begin{align*}
     u=n_{2},~ \tau=\left(\frac{\log\alpha}{\log 2}\right),~ v=a_{1},~ \mu=\frac{\log\left((1+\alpha^{n_{1}-n_{2}})/4\sqrt{2}\right)}{\log 2},
 \end{align*}
with $\left(A,B,w\right)=(6.3,2,({a_{1}-a_{2}})).$
We calculate $\varepsilon>0.00225968.$ Applying Lemma \ref{le:2.2} for $(n_{1}-n_{2})\leq84$, we get $a_{1}-a_{2}\leq 215.$

$\textbf{Step~6}$: We continue case $2.$ We have that $a_{1}-a_{2}\leq 215$ and $n_{1}-n_{2}\leq 84.$ Applying similar steps as in case 1B, we obtain $a_{1}-a_{3}\leq 222.$
We summarize our results obtained so far in the following table.
\begin{center}
\begin{tabular}{l|l|l|l}
Upper bound of & Case 1A & Case 1B & Case 2 \\
\hline
$(a_{1}-a_{2})$ & $214$ & $214$ & $215$\\
$(a_{1}-a_{3})$ &
$218$ & $222$ & $222$ \\
$(n_{1}-n_{2})$ & $87$ &
$86$ & $84$ 
\end{tabular}
\end{center}

$\textbf{Step~7}$: Now, under the assumption that $n_{1}-n_{2}\leq 87, a_{1}-a_{2}\leq 215, a_{1}-a_{3}\leq 222 $, put
 \begin{align*}
    \Lambda_{3}=-n_{1}\log\alpha+a_{1}\log 2+\log\frac{4\sqrt{2}(1+2^{a_{2}-a_{1}}+2^{a_{3}-a_{1}})}{(1+\alpha^{n_{2}-n_{1}})}.
\end{align*}
The inequality \eqref{en:B21} can be written as
\begin{align*}
\left\lvert e^{\Lambda_{3}}-1\right\rvert=\lvert \Gamma_{3}\rvert<0.6\alpha^{-n_{1}}.
\end{align*}
 which implies that
 \begin{align*}
    \left\lvert n_{1}\log\alpha-a_{1}\log 2+\log\frac{(1+\alpha^{n_{2}-n_{1}})}{4\sqrt{2}(1+2^{a_{2}-a_{1}}+2^{a_{3}-a_{1}})}\right\rvert<1.2\alpha^{-n_{1}}.
 \end{align*}
 Dividing both sides by $\log 2$ gives
 \begin{align*}
     \left\lvert n_{1}\left(\frac{\log\alpha}{\log 2}\right)-a_{1}+\frac{\log\left((1+\alpha^{n_{2}-n_{1}})/(4\sqrt{2}(1+2^{a_{2}-a_{1}}+2^{a_{3}-a_{1}}))\right)}{\log 2}\right\rvert&<\frac{1.2}{\log2}\alpha^{-n_{1}}\\
     &<1.7\alpha^{-n_{1}}.
 \end{align*}
Let
 \begin{align*}
     u=n_{1},~ \tau=\left(\frac{\log\alpha}{\log 2}\right),~ v=a_{1},~ \mu=\frac{\log\left((1+\alpha^{n_{2}-n_{1}})/(4\sqrt{2}(1+2^{a_{2}-a_{1}}+2^{a_{3}-a_{1}}))\right)}{\log 2},
 \end{align*}
with $\left(A,B,w\right)=(1.7,\alpha,n_{1}).$
With the same $M$, we find $\varepsilon>0.00001.$ Applying Lemma \ref{le:2.2} for $n_{1}-n_{2}\leq 87, a_{1}-a_{2}\leq 215$  and $a_{1}-a_{3}\leq 222$, we get $n_{1}\leq 86$, which is a contradiction.
Hence, the theorem is proved.

As a consequence of Theorem \ref{th:1} we obtain the following corollaries.
\begin{theorem}
All non-negative integer solutions $(n_{1},n_{2},a_{1},a_{2})$ of the equation
\begin{equation*}
 B_{n_{1}}+B_{n_{2}}=2^{a_{1}}+2^{a_{2}},
\end{equation*} with $n_{1}\geq n_{2}\geq0$ and $a_{1}\geq a_{2} \geq 0$ are given by
\begin{align*}
(n_{1},n_{2},a_{1},a_{2})\in
\{(1,1,0,0), (2,0,2,1), (2,2,3,2),
(3,1,5,2)\}.
\end{align*}
\end{theorem}
\begin{theorem}
All non-negative integer solutions $(n_{1},n_{2},a_{1})$ of the equation
\begin{equation*}
 B_{n_{1}}+B_{n_{2}}=2^{a_{1}},
\end{equation*} with $n_{1}\geq n_{2}\geq  0$ and $a_{1}\geq 0$ are given by
\begin{align*}
(n_{1},n_{2},a_{1})\in
\{(1,1,0), (1,1,1)\}.
\end{align*}
\end{theorem}




\EditInfo{August 21, 2021}{December 04, 2022}{Attila Bérczes}

\end{paper}